\documentclass[11pt]{article}

\usepackage{latexsym}
\usepackage[dvips]{graphicx}
\usepackage{float}
\usepackage{amsmath}
\usepackage{amssymb}
\usepackage{epic}
\usepackage{color}
\usepackage{setspace}
\usepackage{lscape}
\usepackage{arydshln}
\usepackage{fancybox}   

\usepackage{latexsym}
\begin{document}
\bibliographystyle{plain}
\floatplacement{table}{H}
\newtheorem{definition}{Definition}[section]
\newtheorem{lemma}{Lemma}[section]
\newtheorem{theorem}{Theorem}[section]
\newtheorem{corollary}{Corollary}[section]
\newtheorem{proposition}{Proposition}[section]
\newcommand{\sni}{\sum_{i=1}^{n}}
\newcommand{\snj}{\sum_{j=1}^{n}}
\newcommand{\smj}{\sum_{j=1}^{m}}
\newcommand{\sumjm}{\sum_{j=1}^{m}}
\newcommand{\bdis}{\begin{displaymath}}
\newcommand{\edis}{\end{displaymath}}
\newcommand{\beq}{\begin{equation}}
\newcommand{\eeq}{\end{equation}}
\newcommand{\beqn}{\begin{eqnarray}}
\newcommand{\eeqn}{\end{eqnarray}}
\newcommand{\qed}{{\large $\sqcap$ \hskip -0.37cm $\sqcup$}}
\newcommand{\defeq}{\stackrel{\triangle}{=}}
\newcommand{\simleq}{\stackrel{<}{\sim}}
\newcommand{\sep}{\;\;\;\;\;\; ; \;\;\;\;\;\;}
\newcommand{\real}{\mbox{$ I \hskip -4.0pt R $}}
\newcommand{\complex}{\mbox{$ I \hskip -6.8pt C $}}
\newcommand{\integ}{\mbox{$ Z $}}
\newcommand{\realn}{\real ^{n}}
\newcommand{\sqrn}{\sqrt{n}}
\newcommand{\sqrtwo}{\sqrt{2}}
\newcommand{\prf}{{\bf Proof. }}

\newcommand{\onehlf}{\frac{1}{2}}
\newcommand{\thrhlf}{\frac{3}{2}}
\newcommand{\fivhlf}{\frac{5}{2}}
\newcommand{\onethd}{\frac{1}{3}}
\newcommand{\lb}{\left ( }
\newcommand{\lcb}{\left \{ }
\newcommand{\lsb}{\left [ }
\newcommand{\labs}{\left | }
\newcommand{\rb}{\right ) }
\newcommand{\rcb}{\right \} }
\newcommand{\rsb}{\right ] }
\newcommand{\rabs}{\right | }
\newcommand{\lnm}{\left \| }
\newcommand{\rnm}{\right \| }
\newcommand{\lambdab}{\bar{\lambda}}
%
%
\newcommand{\xj}{x_{j}}
\newcommand{\xjb}{\bar{x}_{j}}
\newcommand{\xro}{x_{\resh}}
\newcommand{\xrob}{\bar{x}_{\resh}}
\newcommand{\xsig}{x_{\sigma}}
\newcommand{\xsigb}{\bar{x}_{\sigma}}
\newcommand{\xnmjb}{\bar{x}_{n-j+1}}
\newcommand{\xnmj}{x_{n-j+1}}
\newcommand{\aroj}{a_{\resh j}}
\newcommand{\arojb}{\bar{a}_{\resh j}}
\newcommand{\aroro}{a_{\resh \resh}}
\newcommand{\amuro}{a_{\mu \resh}}
\newcommand{\amumu}{a_{\mu \mu}}
\newcommand{\aii}{a_{ii}}
\newcommand{\aik}{a_{ik}}
\newcommand{\akj}{a_{kj}}
\newcommand{\atwoii}{a^{(2)}_{ii}}
\newcommand{\atwoij}{a^{(2)}_{ij}}
\newcommand{\ajj}{a_{jj}}
\newcommand{\aiib}{\bar{a}_{ii}}
\newcommand{\ajjb}{\bar{a}_{jj}}
\newcommand{\bii}{a_{jj}}
\newcommand{\biib}{\bar{a}_{jj}}
\newcommand{\aij}{a_{i,n-i+1}}
\newcommand{\akl}{a_{j,n-j+1}}
\newcommand{\aijb}{\bar{a}_{i,n-i+1}}
\newcommand{\aklb}{\bar{a}_{j,n-j+1}}
\newcommand{\bij}{a_{n-j+1,j}}
\newcommand{\arorob}{\bar{a}_{\resh \resh}}
\newcommand{\arosig}{a_{\resh \sigma}}
\newcommand{\arosigb}{\bar{a}_{\resh \sigma}}
\newcommand{\sumjrosig}{\sum_{\stackrel{j=1}{j\neq\resh,\sigma}}^{n}}
\newcommand{\summuro}{\sum_{\stackrel{j=1}{j\neq\mu,\resh}}^{n}}
\newcommand{\sumjnoti}{\sum_{\stackrel{j=1}{j\neq i}}^{n}}
\newcommand{\sumlnoti}{\sum_{\stackrel{\ell=1}{\ell \neq i}}^{n}}
\newcommand{\sumknoti}{\sum_{\stackrel{k=1}{k\neq i}}^{n}}
\newcommand{\sumknotij}{\sum_{\stackrel{k=1}{k\neq i,j}}^{n}}
\newcommand{\sumk}{\sum_{k=1}^{n}}
\newcommand{\snl}{\sum_{\ell=1}^{n}}
\newcommand{\sumji}{\sum_{\stackrel{j=1}{j\neq i, n-i+1}}^{n}}
\newcommand{\sumki}{\sum_{\stackrel{k=1}{k\neq i, n-i+1}}^{n}}
\newcommand{\sumkj}{\sum_{\stackrel{k=1}{k\neq j, n-j+1}}^{n}}
\newcommand{\sumjro}{\sum_{\stackrel{j=1}{j\neq\resh}}^{n}}
\newcommand{\rrosig}{R''_{\resh \sigma}}
\newcommand{\rro}{R'_{\resh}}
\newcommand{\gamror}{\Gamma_{\resh}^{R}(A)}
\newcommand{\gamir}{\Gamma_{i}^{R}(A)}
\newcommand{\gamctrr}{\Gamma_{\frac{n+1}{2}}^{R}(A)}
\newcommand{\gamctrc}{\Gamma_{\frac{n+1}{2}}^{C}(A)}
\newcommand{\gamroc}{\Gamma_{\resh}^{C}(A)}
\newcommand{\gamjc}{\Gamma_{j}^{C}(A)}
\newcommand{\lamror}{\Lambda_{\resh}^{R}(A)}
\newcommand{\lamir}{\Lambda_{i}^{R}(A)}
\newcommand{\lamirepsilon}{\Lambda_{i}^{R}(A_{\epsilon})}
\newcommand{\lamnir}{\Lambda_{n-i+1}^{R}(A)}
\newcommand{\lamjr}{\Lambda_{j}^{R}(A)}
\newcommand{\phiij}{\Phi_{ij}^{R}(A)}
\newcommand{\delir}{\Delta_{i}^{R}(A)}
\newcommand{\vir}{V_{i}^{R}(A)}
\newcommand{\pamir}{\Pi_{i}^{R}(A)}
\newcommand{\xir}{\Xi_{i}^{R}(A)}
\newcommand{\lamjc}{\Lambda_{j}^{C}(A)}
\newcommand{\vjc}{V_{j}^{C}(A)}
\newcommand{\pamjc}{\Pi_{j}^{C}(A)}
\newcommand{\xjc}{\Xi_{j}^{C}(A)}
\newcommand{\lamroc}{\Lambda_{\resh}^{C}(A)}
\newcommand{\lamsigr}{\Lambda_{\sigma}^{R}(A)}
\newcommand{\lamsigc}{\Lambda_{\sigma}^{C}(A)}
\newcommand{\psii}{\Psi_{i}^{R}(A)}
\newcommand{\psiq}{\Psi_{q}}
\newcommand{\psiiepsilon}{\Psi_{i}^{R}(A_{\epsilon})}
\newcommand{\psiqepsilon}{\Psi_{q}(A_{\epsilon})}
\newcommand{\psiqc}{\Psi_{q}^{c}}
\newcommand{\psiqcepsilon}{\Psi_{q}^{c}(A_{\epsilon})}

\newcommand{\xmu}{x_{\mu}}
\newcommand{\xmub}{\bar{x}_{\mu}}
\newcommand{\xnu}{x_{\nu}}
\newcommand{\xnub}{\bar{x}_{\nu}}
\newcommand{\amuj}{a_{\mu j}}
\newcommand{\amujb}{\bar{a}_{\mu j}}
\newcommand{\amumub}{\bar{a}_{\mu \mu}}
\newcommand{\amunu}{a_{\mu \nu}}
\newcommand{\amunub}{\bar{a}_{\mu \nu}}
\newcommand{\sumjmunu}{\sum_{\stackrel{j=1}{j\neq\mu,\nu}}}
\newcommand{\rmunu}{R''_{\mu \nu}}
\newcommand{\rmu}{R'_{\mu}}

\newcommand{\Azero}{A_{0}}
\newcommand{\Aone}{A_{1}}
\newcommand{\Atwo}{A_{2}}
\newcommand{\Ath}{A_{3}}
\newcommand{\Afr}{A_{4}}
\newcommand{\Afv}{A_{5}}
\newcommand{\Asx}{A_{6}}
\newcommand{\Anmo}{A_{n-1}}
\newcommand{\Anmt}{A_{n-2}}
\newcommand{\An}{A_{n}}

\newcommand{\azero}{a_{0}}
\newcommand{\aone}{a_{1}}
\newcommand{\atwo}{a_{2}}
\newcommand{\ath}{a_{3}}
\newcommand{\afr}{a_{4}}
\newcommand{\afv}{a_{5}}
\newcommand{\asx}{a_{6}}
\newcommand{\anmo}{a_{n-1}}
\newcommand{\anmt}{a_{n-2}}
\newcommand{\an}{a_{n}}

\newcommand{\Bzero}{B_{0}}
\newcommand{\Bone}{B_{1}}
\newcommand{\Btwo}{B_{2}}
\newcommand{\Bth}{B_{3}}
\newcommand{\Bfr}{B_{4}}
\newcommand{\Bfv}{B_{5}}
\newcommand{\Bsx}{B_{6}}
\newcommand{\Bnmo}{B_{n-1}}
\newcommand{\Bndtwomo}{B_{n/2-1}}
\newcommand{\Bnmt}{B_{n-2}}
\newcommand{\Bn}{B_{n}}

\newcommand{\cii}{c_{ii}}
\newcommand{\cik}{c_{ik}}
\newcommand{\ckj}{c_{kj}}
\newcommand{\ctwoii}{c^{(2)}_{ii}}
\newcommand{\ctwoij}{c^{(2)}_{ij}}
\newcommand{\cjj}{c_{jj}}

\newcommand{\bik}{b_{ik}}
\newcommand{\bkj}{b_{kj}}
\newcommand{\btwoii}{b^{(2)}_{ii}}
\newcommand{\btwoij}{b^{(2)}_{ij}}
\newcommand{\bjj}{b_{jj}}

\newcommand{\abii}{(AB)_{ii}}
\newcommand{\abil}{(AB)_{i\ell}}

\newcommand{\bkl}{b_{k\ell}}
\newcommand{\btwoil}{b^{(2)}_{i\ell}}
\newcommand{\bll}{b_{\ell \ell}}

\newcommand{\matrixspace}{\;\;}
\newcommand{\ellone}{\ell_{1}}
\newcommand{\elltwo}{\ell_{2}}

\newcommand{\phik}{\phi_{k}}
\newcommand{\chik}{\chi_{k}}
\newcommand{\Phik}{\Phi_{k}}
\newcommand{\psik}{\psi_{k}}
\newcommand{\dr}{\beta^{n-k}}
\newcommand{\dn}{\beta^{-k}}
\newcommand{\betatwotilde}{\tilde{\rtwonk}}
\newcommand{\betathtilde}{\tilde{\ronenk}}
\newcommand{\betamink}{\beta^{-k}}

\newcommand{\rphik}{r(\phi_{k})}    
\newcommand{\rrphik}{(r(\phi_{k}))}    
\newcommand{\rchik}{r(\chi_{k})}    
\newcommand{\rrchik}{(r(\chi_{k}))}    
\newcommand{\sphik}{s(\phi_{k})}    
\newcommand{\ssphik}{(s(\phi_{k}))}    
\newcommand{\schik}{s(\chi_{k})}    
\newcommand{\sschik}{(s(\chi_{k}))}    

\newcommand{\rpsik}{r(\psi_{k})}    
\newcommand{\rrpsik}{(r(\psi_{k}))} 
\newcommand{\spsik}{s(\psi_{k})}    
\newcommand{\sspsik}{(s(\psi_{k}))} 
\newcommand{\rpsinmk}{r(\psi_{n-k})}    
\newcommand{\rrpsinmk}{(r(\psi_{n-k}))} 
\newcommand{\spsinmk}{s(\psi_{n-k})}    
\newcommand{\sspsinmk}{(s(\psi_{n-k}))} 

\newcommand{\rphinmk}{r(\phi_{n-k})}    
\newcommand{\rrphinmk}{(r(\phi_{n-k}))}    
\newcommand{\rchinmk}{r(\chi_{n-k})}    
\newcommand{\rrchinmk}{(r(\chi_{n-k}))}    
\newcommand{\sphinmk}{s(\phi_{n-k})}    
\newcommand{\ssphinmk}{(s(\phi_{n-k}))}    
\newcommand{\schinmk}{s(\chi_{n-k})}    
\newcommand{\sschinmk}{(s(\chi_{n-k}))}    
\newcommand{\stilde}{\tilde{s}}

\newcommand{\resh}{\rho}
\newcommand{\snk}{s}
\newcommand{\ronenk}{r_{1}}
\newcommand{\rtwonk}{r_{2}}
\newcommand{\obar}{\widebar{O}}
\newcommand{\Lbar}{\widebar{L}}
\newcommand{\dltaone}{\delta_{1}}
\newcommand{\dltatwo}{\delta_{2}}
\newcommand{\dltatld}{\tilde{\delta}}
\newcommand{\tauone}{\tau_{1}}
\newcommand{\tautwo}{\tau_{2}}
\newcommand{\xb}{\bar{x}}
\newcommand{\qone}{q_{1}}
\newcommand{\qtwo}{q_{2}}
\newcommand{\ub}{\bar{u}}
\newcommand{\cmm}{\complex^{m \times m}}
\newcommand{\opdsk}{\mathcal{O}}
\newcommand{\cldsk}{\widebar{\mathcal{O}}}

\begin{center}
\large
{\bf GENERALIZATION AND VARIATIONS OF PELLET'S THEOREM FOR MATRIX POLYNOMIALS} \\
\vskip 0.5cm
\normalsize
A. Melman \\
Department of Applied Mathematics \\
School of Engineering, Santa Clara University  \\
Santa Clara, CA 95053  \\
e-mail : amelman@scu.edu \\
\vskip 0.5cm
\end{center}

\begin{abstract}
We derive a generalized matrix version of Pellet's theorem, itself based on a 
generalized Rouch\'{e} theorem for matrix-valued functions, to generate upper, lower, and 
internal bounds on the eigenvalues of matrix polynomials. Variations of the theorem are 
suggested to try and overcome situations where Pellet's theorem cannot be applied.
\vskip 0.15cm
{\bf AMS(MOS) subject classification :} 12D10, 15A18, 30C15
\end{abstract}

%
%
%
%

\section{Introduction}           
\label{introduction} 
Polynomial eigenvalue problems have been investigated for quite some time (\cite{GohbergLancasterRodman},
\cite{Lancaster1960}, \cite{Lancaster1966}) and have important applications in a wide range of 
engineering fields such as vibration analysis, acoustics, and fluid mechanics - to name just a few
(\cite{TisseurMeerbergen}).
It is, in general, costly to compute polynomial eigenvalues for large problems, 
but bounds on such eigenvalues are relatively easy to obtain. They provide information on the location 
of eigenvalues that can be used by iterative methods for computing them and are also useful in the 
computation of pseudospectra.

The polynomial eigenvalue problem is to find a nonzero eigenvector $v$, corresponding to an eigenvalue $z$
satisfying $P(z)v=0$, where
\bdis
P(z) = A_{n} z^{n} + A_{n-1} z^{n-1} + \dots + A_{0},
\edis
with $A_{j} \in \complex^{m \times m}$ for $j=0,\dots n$. We will assume throughout that $\det{(P)}$ is not
identically zero.
If $A_{n}$ is singular then $P$ has infinite eigenvalues and if $A_{0}$ is singular then zero is an eigenvalue.
There are $nm$ eigenvalues, including possibly infinite ones. The finite eigenvalues are the solutions of
$\det{(P)}=0$.

To explain the aims of this work, we start with the scalar polynomial with complex coefficients
\bdis
p(z) = a_{n} z^{n} + a_{n-1} z^{n-1} + \dots + a_{0} \; ,
\edis
and $a_{n} a_{0} \neq 0$. If $\zeta$ is a zero of $p$, then for $a_{k} \neq 0$ with $1 \leq k \leq n-1$ we have
\bdis
-a_{k} \zeta^{k} = a_{n} \zeta^{n} + a_{n-1} \zeta^{n-1} + \dots 
+ a_{k+1} \zeta^{k+1} + a_{k-1} \zeta^{k-1} + \dots + a_{0}\; ,
\edis
and an analogous equality for $k=0,n$, which implies that
\begin{eqnarray}
|a_{k}| \, | \zeta|^{k} &   =  & \left | a_{n} \zeta^{n} + a_{n-1} \zeta^{n-1} + \dots 
+ a_{k+1} \zeta^{k+1} + a_{k-1} \zeta^{k-1} + \dots + a_{0} \right |  \nonumber \\
                        & \leq & |a_{n}| \, | \zeta|^{n} + |a_{n-1}| \, |\zeta|^{n-1} + \dots 
+ |a_{k+1}| \, | \zeta|^{k+1} + |a_{k-1}| \,  |\zeta|^{k-1} + \dots + |a_{0}| \; .  \label{introineq1} 
\end{eqnarray}
Applying the equivalent of inequality~(\ref{introineq1}) with $k=n$, we obtain that $|\zeta |$ must satisfy
\beq
\label{Cauchy1}
|a_{n}| \, | \zeta|^{n} - |a_{n-1}| \, |\zeta|^{n-1} -  \dots - |a_{1}| |\zeta | - |a_{0}| \leq  0 \; .  
\nonumber
\eeq
This means that $|\zeta |$ can be no larger than the unique positive root of
\beq
\label{Cauchy2}
| a_{n}| x^{n} - |a_{n-1}| x^{n-1} - \dots -|a_{1}| x - |a_{0}| = 0 \; . 
\nonumber
\eeq
Analogously, one finds that $|\zeta |$ can be no smaller than the unique positive root of
\beq
\label{Cauchy3}
| a_{n}| x^{n} + |a_{n-1}| x^{n-1} + \dots +|a_{1}| x - |a_{0}| = 0 \; . 
\nonumber
\eeq
These results are due to Cauchy (\cite{Cauchy}, \cite[Theorem (27,1), p.122]{Marden}). 

When we apply inequality~(\ref{introineq1}) with $1 \leq k \leq n-1$ and $a_{k} \neq 0$, we obtain 
\beq
\label{Pellet1}
|a_{n}| \, | \zeta|^{n} + |a_{n-1}| \, |\zeta|^{n-1} + \dots 
+ |a_{k+1}| \, | \zeta|^{k+1} - |a_{k}| \, | \zeta|^{k} + |a_{k-1}| \,  |\zeta|^{k-1} + \dots + |a_{0}|  \geq 0 \; .
\eeq
Now, the equation
\beq
\label{Pellet2}
|a_{n}| x^{n} + |a_{n-1}| x^{n-1} + \dots + |a_{k+1}| x^{k+1} - |a_{k}| x^{k} + |a_{k-1}| x^{k-1} 
+ \dots + |a_{0}| = 0 
\nonumber
\eeq
has, by Descartes' rule, either two or no positive roots. 
If it has no positive roots or if the two positive roots coincide,
then inequality~(\ref{Pellet1}) is always satisfied and provides no information on $|\zeta |$. However, 
if it has two distinct positive roots, say $x_{1}$ and $x_{2}$ 
with $x_{1} < x_{2}$, then inequality~(\ref{Pellet1}) implies that either $|\zeta | \leq x_{1}$ or 
$|\zeta | \geq x_{2}$ and that no zero of $p$ has a modulus in $(x_{1},x_{2})$. As a straightforward consequence  
of Rouch\'{e}'s theorem, 
Pellet's theorem makes this result more precise. 
It is stated as follows.
\begin{theorem}
({\bf Pellet}, \cite{Pellet}, \cite[Th.(28,1), p.128]{Marden})
\label{Pellet}
Given the polynomial $p(z)=z^{n} + a_{n-1} z^{n-1} + \dots + a_{1} z + a_{0}$ with complex coefficients, 
$a_{0} a_{k} \neq 0$, and $n \geq 3$, let $1 \leq k \leq n-1$, and let the polynomial 
\bdis
f_{k}(x)=x^{n} + |a_{n-1}| x^{n-1} + \dots + |a_{k+1}| x^{k+1} - |a_{k}| x^{k} + |a_{k-1}|x^{k-1} + \dots  + |a_{0}| 
\edis
have two distinct positive roots $x_{1}$ and $x_{2}$ with $x_{1} < x_{2}$.
Then $p$ has exactly $k$ zeros in or on the circle $|z| = x_{1}$ and no zeros in the annular ring 
$x_{1} < |z| < x_{2}$.
\end{theorem}
Cauchy's results can be considered as a special limit case of Pellet's theorem.

We have two main goals. The first is to generalize Pellet's theorem (and therefore also Cauchy's results) to 
matrix polynomials, with matrix norms replacing absolute values. Such a generalized theorem for the spectral
norm (2-norm) was recently derived in \cite{BiniNoferiniSharify}, where it was used to determine initial 
approximations for iterative methods like the Ehrlich-Aberth method (\cite{Aberth}, \cite{BiniNumAlg}, \cite{Ehrlich}).
Using a different proof, we obtain a generalization valid for any subordinate norm, not just the spectral norm.

Our second goal is to find a way to, at least sometimes, overcome situations where $f_{k}$ in Theorem~\ref{Pellet} 
(or an analogous polynomial for the generalized version of that theorem) does not have positive roots for a 
particular value of $k$. To do this, we derive a variation of the generalized Pellet theorem that relies on the 
existence of positive roots for a polynomial of roughly half the degree of $f_{k}$.

Only polynomials (or matrix polynomials) with very special coefficients would have more than a few values of $k$ for 
which Pellet's theorem can be applied and even when there exists a significant gap between groups of zeros 
(or eigenvalues), it frequently occurs that the theorem is unable to detect it. 
Unfortunately, there exist no results for Pellet's theorem that would allow one to predict if the function $f_{k}$
has positive roots or not, and the same problem naturally carries over to our variation of that theorem. This makes 
it impossible to predict if and when our result improves over Pellet's theorem. On the other hand, it does provide an 
alternative that currently does not exist in cases where Pellet's theorem is not applicable and we will present extensive 
numerical examples illustrating its usefulness.

We also suggest an idea applicable to scalar polynomials, whereby the zeros of a polynomial are computed as
the eigenvalues of an appropriate matrix polynomial. Perhaps paradoxically, this can lead (like the variation 
mentioned above) to an improvement of the original (scalar) version of Pellet's theorem by using 
its generalized (matrix) version. At the same time, this approach can also improve the upper and lower bounds 
resulting from Cauchy's result. Its potential will be demonstrated by numerical 
experiments for a restricted class of polynomials.

We will not dwell on the numerical problem of finding the positive roots of functions like $f_{k}$ in 
Theorem~\ref{Pellet}, if they exist. An efficient way to compute them can be found in \cite{MelmanPelletImpl},
while bounds on such roots were derived in \cite{BiniNoferiniSharify}. The number of arithmetic operations that
this requires is typically dwarfed by the much more costly computation of the eigenvalues of a matrix polynomial.
Moreover, this work focuses on Pellet's theorem \emph{if and when} it is used.
It does not focus on whether or not it should be used in the first place, as this depends very
much on external factors.
  
The organization of the paper is as follows. In Section~2 we review a few preliminaries that will be needed
in the later sections. In Section~3 we derive the generalized Pellet theorem, in Section~4 we propose
a variation of that theorem, and in Section~5 we present the aforementioned idea for scalar polynomials.


%
%
%
%
\section{Preliminaries}
\label{preliminaries}

Since it is mentioned several times, we begin by stating Rouch\'{e}'s well-known theorem.
%
%
\begin{theorem}
\label{Rouchetheorem}
({\bf Rouch\'{e}}, \cite{Rouche}, \cite[Theorem 1.6]{Lang})
Let $f$ and $g$ be analytic in the interior of a simple closed curve $\Gamma$ and continuous on $\Gamma$, and let
$|g(z)| < |f(z)|$ for all $z \in \Gamma$. Then $f+g$ and $f$ have the same number of zeros in the interior of the curve.
\end{theorem}
All matrix norms throughout this work are assumed to be vector-induced (or \emph{subordinate}). 
\nopagebreak The main matrix norms we will use are the $1$-Norm, $\infty$-Norm, and the $2$-Norm 
(or \emph{spectral norm}), defined (\cite[p. 294-295]{HJ}) for $A \in \complex^{n \times n}$
with elements $a_{ij}$ by
\begin{eqnarray*}
|| A ||_{1}      & = & \max_{1 \leq j \leq n} \sum_{i=1}^{n} |a_{ij}| = || A^{*} ||_{\infty}  \; , \\
|| A ||_{\infty} & = & \max_{1 \leq i \leq n} \sum_{j=1}^{n} |a_{ij}| = || A^{*} ||_{1} \; , \\
|| A ||_{2} & = & \max \left \{ \sqrt{\lambda} \; : \; \lambda \; \text{is an eigenvalue of} \;  A^{*}A \right \}
 = || A^{*} ||_{2}\; . 
\end{eqnarray*}
The zeros of a monic polynomial $p(z) = z^{n} + \anmo z^{n-1} + \dots + \azero$ with complex coefficients
are the eigenvalues of its $n \times n$ companion matrix
\bdis
C(p) =
\begin{pmatrix}
0 &        &       &   & -\azero     \\
1 &        &       &   & -\aone     \\
  & \ddots &       &   & \vdots \\
  &        &       & 1 & -\anmo    \\
\end{pmatrix}
\; . 
\edis
Likewise, the eigenvalues of the monic matrix polynomial 
\bdis
P(z) = I z^{n} + \Anmo z^{n-1} + \dots + \Aone z + \Azero \; , 
\edis
where $A_{j} \in \complex^{m \times m}$ for $j=0,\dots,n-1$, are given by the eigenvalues
of the $nm \times nm$ block companion matrix
\bdis
C(P) =
\begin{pmatrix}
0 &        &       &   & -\Azero     \\
I &        &       &   & -\Aone     \\
  & \ddots &       &   & \vdots \\
  &        &       & I & -\Anmo    \\
\end{pmatrix}
\; , 
\edis
with $I$ the $m \times m$ identity matrix. Since the size of $I$ will usually be clear from the context, 
it will be omitted from the notation.

The reciprocal monic polynomial of $p(z) = z^{n} + \anmo z^{n-1} + \dots + \azero$ with $\azero \neq 0$ is defined by
$p_{r}(z) = z^{n}p(1/z)/a_{0}$. Its zeros are the reciprocals of the zeros of $p$.
Likewise, the reciprocal matrix polynomial of $P(z) = I z^{n} + \Anmo z^{n-1} + \dots + \Aone z + \Azero$ with
$A_{j} \in \complex^{m \times m}$ for $j=0,\dots,n-1$ and $A_{0}$ nonsingular is defined by 
$P_{r}(z) = z^{n}A_{0}^{-1}P(1/z)$. Its eigenvalues are the reciprocals of the eigenvalues of $P$.

%
%
%
%
\section{Generalized Rouch\'{e} and Pellet theorems}
\label{genRouchePellet}

The standard proof of Pellet's theorem uses Rouch\'{e}'s theorem for analytical functions, which suggests that a 
generalized version of Rouch\'{e}'s theorem for analytical matrix-valued functions will be needed. Such a theorem
can be derived from the generalized Rouch\'{e} theorem for bounded linear operators in 
\cite[Theorem 9.2, p. 206]{GohGoldKaas} (see also \cite{GohbergSigal}).
To do this, we need a few definitions for which
we have used similar notation and style as in \cite[Chapter XI]{GohGoldKaas}.
We denote by $\mathcal{L}(X,Y)$ the space of all bounded linear operators from $X$ to $Y$, where $X$ and $Y$ are 
complex Banach spaces. When $X=Y$, this space is denoted by $\mathcal{L}(X)$. 
A bounded linear operator $S : X \rightarrow Y$, where $X$ and $Y$ are complex Banach spaces is called a 
\emph{Fredholm operator} if its range Im $S$ is closed and if dim(Ker $S$) and dim($Y/\text{Im}S$) are finite. 
If $X$ and $Y$ are finite dimensional, then $S$ is always Fredholm. 

For an open connected subset $\Omega$ of $\complex$, let $Q : \Omega \rightarrow \mathcal{L}(X)$ 
be a bounded operator function that is analytic on $\Omega$ and let $\Gamma$ be a simple closed curve in $\Omega$ 
such that its inner domain is a subset of $\Omega$. Then $Q$ is said to be \emph{normal with respect to} 
$\Gamma$ if $Q(z)$ is invertible for all $z \in \Gamma$ and $Q(z)$ is Fredholm on the inner domain of~$\Gamma$. 
Finally, we recall that a Hilbert space is separable if and only if it has a countable orthonormal basis.

The following theorem generalizes Rouch\'{e}'s theorem for operator-valued functions. We have stated it as it
appears in \cite[p. 206]{GohGoldKaas}, but a similar version was already published in \cite{GohbergSigal}.

%
%
\begin{theorem}
\label{genRoucheGohberg}    
{\bf (Generalized Rouch\'{e} theorem for operators.)}
(\cite{GohbergSigal}, \cite{GohGoldKaas})
Let $W,S:\Omega \rightarrow \mathcal{L}(H)$ be analytic operator functions, where $\Omega$ is an open
connected subset of $\complex$ and $H$ is a separable Hilbert space, and assume that 
$W$ is normal with respect to the simple closed curve $\Gamma \subseteq \Omega$. 

If $||W(z)^{-1}S(z)|| < 1$ for all $z \in \Gamma$,
then $W+S$ is also normal with respect to $\Gamma$ and $W+S$ and $W$ have the same number of finite eigenvalues
inside~$\Gamma$, counting multiplicities.
\end{theorem}
The norm used in the previous theorem can be any norm, induced by a norm on $H$.
We obtain the following generalized Rouch\'{e} theorem for analytic matrix-valued functions as an immediate consequence
of Theorem~\ref{genRoucheGohberg}. We note that matrix polynomials are special cases of such functions.
%
%
\begin{theorem}
\label{genRoucheMatrices}           
{\bf (Generalized Rouch\'{e} theorem for matrices.)}
Let $A,B:\Omega \rightarrow \cmm$ be analytic matrix-valued functions, where $\Omega$ is an open
connected subset of $\complex$ 
and assume that $A(z)$ is nonsingular for all $z$ on the simple closed curve $\Gamma \subseteq \Omega$. 

If $||A(z)^{-1}B(z)|| < 1$ for all $z \in \Gamma$,
then $\det(A+B)$ and $\det(A)$ have the same number of zeros inside $\Gamma$, counting multiplicities.
\end{theorem}
\prf
The proof follows directly from Theorem~\ref{genRoucheGohberg}. In this case, $H~\equiv~\complex^{n}$, a finite dimensional
Hilbert space which is (trivially) 
separable and $\mathcal{L}~(H)~\equiv~\cmm$. Furthermore, since $A(z)$ is a matrix, it is a bounded linear 
operator that is normal with respect to $\Gamma$ because $A(z)$ is invertible for all $z \in \Gamma$ and,
since $\complex^{n}$ is finite dimensional, $A(z)$ is (trivially) a Fredholm operator for any $z$ in the interior 
of $\Gamma$.
By Theorem~\ref{genRoucheGohberg}, $A+B$ is then also normal with respect to $\Gamma$ and has the same number of finite 
eigenvalues inside $\Gamma$ as $A$, which means that $\det(A+B)$ and $\det(A)$ have the same number of zeros 
inside $\Gamma$, counting multiplicities. \qed
$ $ \newline
{\bf Remarks.}
{\bf (1)} The norm in Theorem~\ref{genRoucheMatrices} can be any matrix norm, induced by a norm on $\complex^{m}$.
\newline {\bf (2)} Since $||A^{-1}(z)B(z)|| \leq ||A^{-1}(z)|| \cdot ||B(z)||$, it is sufficient that
$||B(z)|| < ||A^{-1}(z)||^{-1}$ for the condition $||A^{-1}(z)B(z) || < 1$ to be satisfied.

Finally, Theorem~\ref{genRoucheMatrices} leads to the following generalization of Pellet's theorem (Theorem~\ref{Pellet}).

%
%
\begin{theorem}
\label{genPellet}             
{\bf (Generalized Pellet theorem.)}
Let 
\bdis  
\label{Pelleteq}
P(z) = A_{n}z^{n} + A_{n-1}z^{n-1} + \dots + A_{1}z + A_{0}  \nonumber 
\edis
be a matrix polynomial with $n \geq 2$, $A_{j} \in \cmm$ for $j=0,\dots,n$, and $A_{0} \neq 0$.      
Let $A_{k}$ be invertible for some $k$ with $1 \leq k \leq n-1$, and let the polynomial 
\bdis
f_{k}(x)=||A_{n}||x^{n} + ||A_{n-1}||x^{n-1} + \dots + ||A_{k+1}||x^{k+1} - ||A^{-1}_{k}||^{-1}x^{k} 
+ ||A_{k-1}||x^{k-1} + \dots + ||A_{1}||x + ||A_{0}||   
\edis
have two distinct positive roots $x_{1}$ and $x_{2}$ with $x_{1} < x_{2}$.
Then $\det(P)$ has exactly $km$ zeros in or on the disk $|z|=x_{1}$ and no zeros in the annular ring 
$x_{1} < |z| < x_{2}$.
\end{theorem}
\prf
By Theorem~\ref{genRoucheMatrices}, if
\beq        
\label{Pelletineq1}
\left |\left |\left (A_{k}z^{k} \right )^{-1} \left ( A_{n} z^{n} + A_{n-1}z^{n-1} + \dots 
+ A_{k+1} z^{k+1} + A_{k-1} z^{k-1} + \dots +  A_{1}z + A_{0} \right ) \right | \right | < 1
\eeq    
for $|z|=x$, then $P$ and $A_{k}z^{k}$ both have the same number of eigenvalues in the open disk $|z| < x$,
namely, $km$. 
From Remark (2) following the proof of Theorem~\ref{genRoucheMatrices}, a sufficient condition for 
inequality~(\ref{Pelletineq1}) to be satisfied is 
\beq
\label{Pelletineq2}
\left | \left | A_{n} z^{n} +  A_{n-1} z^{n-1} + \dots + A_{k+1} z^{k+1} + A_{k-1} z^{k-1} + \dots + A_{1} z + A_{0}
\right | \right | < \left | \left | \left ( A_{k}z^{k} \right )^{-1} \right | \right |^{-1} \; . 
\eeq
Since 
\begin{eqnarray*}
& & \hskip -1.0cm \left |\left | 
A_{n} z^{n} +  A_{n-1} z^{n-1} + \dots + A_{k+1} z^{k+1} + A_{k-1} z^{k-1} + \dots + A_{1} z + A_{0}
\right |\right |   \\
& & \hskip -0.50cm \leq 
||A_{n}|| \, |z|^{n} + ||A_{n-1}|| \, |z|^{n-1} + \dots + ||A_{k+1}|| \, |z|^{k+1} 
+ ||A_{k-1}|| \, |z|^{k-1} + \dots + ||A_{1}|| \, |z| + ||A_{0}|| \; ,
\end{eqnarray*}
inequality~(\ref{Pelletineq2}) will be satisfied for $|z|=x$ whenever $x$ is such that
\beq
\label{Pelletineq3}
||A_{n}||x^{n} + ||A_{n-1}||x^{n-1} + \dots + ||A_{k+1}||x^{k+1} 
+ ||A_{k-1}||x^{k-1} + \dots + ||A_{1}||x + ||A_{0}||   
< ||A^{-1}_{k}||^{-1}x^{k} \; , 
\eeq
i.e., when $f_{k}(x) < 0$, where $f_{k}$ was defined in the statement of the theorem.
By Descartes' rule, $f_{k}$ has either two or no positive roots. Since it was assumed that $f_{k}$ has two positive 
roots $x_{1}$ and $x_{2}$ and since $f_{k}(0) > 0$, inequality~(\ref{Pelletineq3}) will be satisfied for any $x$ 
such that $x_{1} < x < x_{2}$. This concludes the proof. \qed
\newline {\bf Remarks.}
{\bf (1)} Applying Theorem~\ref{genPellet} to the matrix polynomial $A_{k}^{-1}P$ yields better values for
$x_{1}$ and $x_{2}$ since $||A_{k}^{-1}A_{j}|| \leq ||A_{k}^{-1}|| \, ||A_{j}||$. The disadvantage is that all 
the coefficients need to be multiplied by $A_{k}^{-1}$, which could be costly.
\newline {\bf (2)} A version of Theorem~\ref{genPellet} (applied to $A_{k}^{-1}P$) was also obtained 
in \cite{BiniNoferiniSharify} by using a different generalization of 
Rouch\'{e}'s theorem from \cite{MondenArimoto} and \cite{VM}, 
which limits that result to the spectral norm. Because it is based on Theorem~\ref{genRoucheMatrices}, 
Theorem~\ref{genPellet} holds for any matrix norm induced by a vector norm. This becomes more important
as the size of the coefficient matrices increases because the spectral norm can be costly to evaluate,
whereas, e.g., the $1$-Norm and $\infty$-Norm can be computed cheaply. Of course, what constitutes a costly
computation is a matter of opinion and there are undoubtedly cases where using the spectral norm is appropriate.
\newline {\bf (3)} Whether the polynomial $f_{k}$ will have two positive roots depends on the relative sizes
of its coefficients. Generally speaking, this will happen when the coefficient of $x^{k}$ is large
compared to most of the other coefficients, but there exist no results that predict how large it needs to be.
As in the scalar case, if $f_{k}$ does not have positive zeros, 
then Theorem~\ref{genPellet} provides no information about a possible gap between groups of eigenvalues. 

The following theorem, which is a generalization of Cauchy's result, can be considered 
as a special limit case of the generalized Pellet theorem just as in the scalar case. Its proof is entirely analogous
and will therefore be omitted.

%
%
\begin{theorem}
\label{genCauchy}             
{\bf (Generalized Cauchy theorem.)}
All eigenvalues of the matrix polynomial
\bdis  
\label{Cauchyeq}
P(z) = A_{n}z^{n} + A_{n-1}z^{n-1} + \dots + A_{1}z + A_{0} \; ,
\edis
where $A_{j} \in \cmm$, for $j=0,\dots,n$, lie in $|z| \leq R$ when $A_{n}$ is nonsingular, and lie
in $|z| \geq r$ when $A_{0}$ is nonsingular,
where $R$ and $r$ are the unique positive roots of
\bdis
||A^{-1}_{n}||^{-1}x^{n} - ||A_{n-1}||x^{n-1} - \dots - ||A_{1}||x - ||A_{0}|| = 0 
\edis
and
\bdis
||A_{n}||x^{n} + ||A_{n-1}||x^{n-1} + \dots + ||A_{1}||x - ||A^{-1}_{0}||^{-1} = 0 \; ,
\edis
respectively.
\end{theorem}
This theorem was also derived in a different way as Lemma 3.1 in \cite{HighamTisseur}.
Again, it can be applied to $A_{n}^{-1}P$ and $A_{0}^{-1}P$ to yield tighter bounds,
but at the expense of adding $n$ multiplications of two $m \times m$ matrices.

%
%
%
%
\section{Variation of Pellet's theorem}
\label{variation}      


The following theorem is a variation of Pellet's theorem, obtained by squaring the block companion matrix and 
repartitioning the result, which leads to additional Cauchy-type and Pellet-type bounds.

\begin{theorem}
\label{squarematrix}
The squares of the eigenvalues of the monic matrix polynomial
\bdis
P(z) = I z^{n} + \Anmo z^{n-1} + \dots + \Aone z + \Azero \; , 
\edis
where $n$ is a positive even integer and $A_{j} \in \complex^{m \times m}$ for $j=0,\dots,n-1$, are given by the eigenvalues
of the monic matrix polynomial
\bdis
Q(z) = I z^{n/2} + \Bndtwomo z^{n/2-1} + \dots + \Bone z + \Bzero  \; ,
\edis
where $B_{j} \in \complex^{2m \times 2m}$ for $j=0,\dots,n/2-1$, with 
\bdis
B_{0} =
\begin{pmatrix}
A_{0} & -A_{0}A_{n-1} \\
A_{1} & -A_{1}A_{n-1}+A_{0} \\
\end{pmatrix}
\;\; ; \;\;
B_{j} =
\begin{pmatrix}
A_{2j} & -A_{2j}A_{n-1}+A_{2j-1} \\
A_{2j+1} & -A_{2j+1}A_{n-1}+A_{2j} \\
\end{pmatrix}
\;,\; j=1,2, \dots, n/2-1 \; .
\edis
Moreover, 
let $B_{k}$ be invertible for some $k$ with $1 \leq k \leq n/2-1$, and let the polynomial 
\bdis
g_{k}(y)=y^{n/2} + ||B_{n/2-1}||y^{n/2-1} + \dots + ||B_{k+1}||y^{k+1} - ||B^{-1}_{k}||^{-1}y^{k} 
+ ||B_{k-1}||y^{k-1} + \dots + ||B_{1}||y + ||B_{0}||   
\edis
have two distinct positive roots $y_{1}$ and $y_{2}$ with $y_{1} < y_{2}$.
Then $P$ has exactly $2km$ eigenvalues in or on the disk $|z|=\sqrt{y_{1}}$ and no eigenvalues in the annular ring 
$\sqrt{y_{1}} < |z| < \sqrt{y_{2}}$. 

In addition, all eigenvalues of $P$ lie in $|z| \leq \sqrt{\rho}$ and, when $B_{0}$
is nonsingular, in $z \geq \sqrt{\tau}$, where $\rho$ and $\tau$ are the unique positive roots of
\bdis
y^{n/2} - ||B_{n/2-1}||y^{n/2-1} - \dots - ||B_{1}||y - ||B_{0}|| = 0 
\edis
and
\bdis
y^{n/2} + ||B_{n/2-1}||y^{n/2-1} + \dots + ||B_{1}||y - ||B^{-1}_{0}||^{-1} = 0 \; ,
\edis
respectively.
\end{theorem}  
\prf
The companion matrix $C(P)$ of $P$ and its square $C^{2}(P)$, whose eigenvalues are the squares of those of $P$,
are given by
\bdis
C(P) =
\begin{pmatrix}
0 &        &       &   & -\Azero     \\
I &        &       &   & -\Aone     \\
  & \ddots &       &   & \vdots \\
  &        &       & I & -\Anmo    \\
\end{pmatrix}
\;\;\; \text{and} \;\;\;
C^{2}(P) = 
\begin{pmatrix}
0 &         &    & -\Azero      & \Azero \Anmo     \\
0 &         &    & -\Aone  & \Aone \Anmo - \Azero      \\
I &          &    & -\Atwo  & \Atwo \Anmo - \Aone \\
  & \ddots  &    &  \vdots    & \vdots   \\
  &         &  I & -\Anmo & \Anmo^{2}-\Anmt\\
\end{pmatrix} \; .
\edis
Since $n$ is even, $C^{2}(P)$ can be repartitioned into $n/2$ blocks of size $2m \times 2m$ as follows:
\bdis
C^{2}(P) = 
\left (
\;
\begin{matrix}
\boxed{\begin{matrix}
0 &  0 \\
0 &  0 \\
\end{matrix}}
&
\begin{matrix}
  &    \\
  &    \\
\end{matrix}
&
\begin{matrix}
  &    \\
  &    \\
\end{matrix}
&
\boxed{\begin{matrix}
-\Azero & \Azero\Anmo       \\
-\Aone  & \Aone \Anmo  - \Azero  \\
\end{matrix}}
\\
%
%
%
& & & \\
%
%
\boxed{\begin{matrix}
I &  0  \\    
0 &  I  \\   
\end{matrix}}
&
\begin{matrix}
  &   \\
  &   \\
\end{matrix}
&
\begin{matrix}
 &  \\
 &  \\
\end{matrix}
&
\boxed{\begin{matrix}
-\Atwo  & \Atwo \Anmo - \Aone  \\
-\Ath   & \Ath \Anmo - \Atwo  \\
\end{matrix}}
\\
%
%
  & \ddots &         & \vdots  \\
%
%
\begin{matrix}
 & \\
 & \\
\end{matrix}
&
\begin{matrix}
 & \\
 & \\
\end{matrix}
&
\boxed{\begin{matrix}
I & 0 \\
0 & I \\
\end{matrix}}
&
\boxed{\begin{matrix}
-\Anmt  & \Anmt \Anmo - A_{n-3} \\
-\Anmo   & \Anmo^{2}-\Anmt  \\
\end{matrix}}
\\
\end{matrix}
\;
\right )
\; .
\edis

In other words, $C^{2}(P)$ can be written as
\beq
\label{blockcompanionCPsq}
C^{2}(P) =
\begin{pmatrix}
0 &        &       &   & -\Bzero     \\
I &        &       &   & -\Bone     \\
  & \ddots &       &   & \vdots \\
  &        &       & I & -B_{n/2-1}    \\
\end{pmatrix}
\; ,
\nonumber
\eeq
where
\bdis
B_{0} =
\begin{pmatrix}
A_{0} & -A_{0}A_{n-1} \\
A_{1} & -A_{1}A_{n-1}+A_{0} \\
\end{pmatrix}
\;\; ; \;\;
B_{j} =
\begin{pmatrix}
A_{2j} & -A_{2j}A_{n-1}+A_{2j-1} \\
A_{2j+1} & -A_{2j+1}A_{n-1}+A_{2j} \\
\end{pmatrix}
\;,\; j=1,2,\dots, n/2-1 \; ,
\edis
and where $I$ now stands for the $2m \times 2m$ identity matrix. The expression for
$C^{2}(P)$ in~(\ref{blockcompanionCPsq}) is the block companion matrix $C(Q)$ of the matrix polynomial
\beq    
\label{qequation}
Q(z) = I z^{n/2} + \Bndtwomo z^{n/2-1} + \dots + \Bone z + \Bzero  \; ,
\eeq    
whose eigenvalues, being the same as those of $C^{2}(P)$, are the squares of the eigenvalues of $P$. 
The remainder of the proof follows directly from applying Theorem~\ref{genPellet} and Theorem~\ref{genCauchy} to the matrix 
polynomial $Q$ and the fact that its eigenvalues are the squares of those of $P$. \qed

Theorem~\ref{squarematrix} can be used to potentially improve over Theorem~\ref{genPellet} for a given
matrix polynomial~$P$. This may take the form of improved upper and/or
lower bounds, but it may also be the case that 
a gap between groups of eigenvalues can now be computed that would otherwise not have been detected,
which is often more important.
Furthermore, the degree of $Q$ is only half that of $P$ which has a generally beneficial effect
on computations involving it.
Of course, the cost of squaring $C(P)$, while still orders of magnitude lower than that of actually computing the
eigenvalues, may not be negligable, depending on the size of $m$.

Additional bounds can be derived by applying Theorem~\ref{squarematrix} to the reciprocal matrix polynomial 
$P_{r}$, which yields
a matrix polynomial that we designate by $Q_{R}$. This matrix polynomial is, in general, different 
from $Q_{r}$ which is the reciprocal polynomial of $Q$, defined in~(\ref{qequation}).  
The reciprocal of the Cauchy upper bound for $P_{r}$ or $Q_{r}$ is equal to the Cauchy lower bound 
for $P$ or $Q$, respectively. However, the reciprocal of the Cauchy upper bound for $Q_{R}$, while being a 
lower bound on the moduli of the eigenvalues of $P$, is, in general, different from the Cauchy lower bound for $Q$.
An analogous situation exists for the reciprocal of the Cauchy lower bound for $Q_{R}$. 
We point out that $P$, $P_{r}$, $Q$, $Q_{r}$, and $Q_{R}$ here are \emph{monic} polynomials.

When $n$ is odd we consider $zP(z)$ instead of $P(z)$, which simply has $m$ extra zero eigenvalues. 
In such a case, $A_{0}$ is the zero matrix, which makes $B_{0}$ in Theorem~\ref{squarematrix}
singular so that $Q$ cannot be used to
obtain a lower bound. To remedy this, one can use $zP_{r}(z)$: the reciprocal of the upper bound 
on its largest eigenvalue provides the desired lower bound on the smallest eigenvalue of $P$.

Many different additional bounds can be derived by considering $MP$ or $NQ$ instead of $P$ or $Q$ for any 
nonsingular matrices $M$ and $N$, or by using different scalings of the companion matrices. It would therefore 
be impractical to compare all possible variants. It is also difficult if not impossible to predict which 
particular bound will outperform other bounds for any given situation. 
Instead, we will present several numerical examples to illustrate that our bounds can be valuable complements to 
existing ones.

The bounds' complexity for a particular value of $k$ is linear in the degree $n$ and further depends on the coefficient
matrices' properties. If, for instance, the coefficients exhibit structure, such as symmetry or sparsity, 
then the various matrix manipulations involved in their computation generally require fewer operations. 
The choice of bounds therefore ultimately depends on the situation, which is outside the scope of our discussion.

We observe that Theorem~\ref{squarematrix} leads to the slightly counterintuitive situation where the 
generalized matrix version of the theorem is applied to a scalar polynomial.

Theorem~\ref{squarematrix} always provides upper and lower bounds 
but it does have limitations since the $2m \times 2m$ matrix coefficients of $Q$ 
are twice the size of those of $P$. This leads to the computation of gaps between groups of eigenvalues that contain an 
even multiple of $m$, instead of just a multiple of $m$. 

Finally, we remark that higher powers of $C(P)$ could similarly be used, although the computational cost involved
may render the resulting bounds inefficient.

We now consider three examples to compare Theorem~\ref{genPellet} and Theorem~\ref{genCauchy} with 
Theorem~\ref{squarematrix}.
In Example 1 we compare the upper and lower bounds resulting from Theorem~\ref{genCauchy} with those from 
Theorem~\ref{squarematrix} for both scalar and matrix polynomials. In Example 2
we compare the number of gaps between groups of eigenvalues detected by Theorem~\ref{genPellet} 
with the number of gaps detected by Theorem~\ref{squarematrix} for matrix polynomials and in Example 3 we 
do the same for scalar polynomials. We also compared the gap width, although the ability to detect a gap
is usually more important than improving the gap width, and it is our first priority.
These examples focus somewhat more on Pellet bounds than on Cauchy bounds since they
are less frequently mentioned in the literature.
Although $n$ is an even number in the examples, this is not essential and
similar conclusions can be drawn for odd numbers.

$ $ \newline \underline{\bf Example 1}
\newline
\newline For this example we have generated matrix polynomials of degree $n=10$ with $m=1,2,10,25$, where
$m=1$ corresponds to a scalar polynomial. They take the form 
\bdis
P(z) = Iz^{10} + A_{9}z^{9} + \dots + A_{0} \; ,
\edis
where each of the $m \times m$ matrices $A_{j}$, for $j=0,\dots ,9$, has complex elements whose real and imaginary 
parts are uniformly distributed in $[-1,1]$, multiplied by a random number which is uniformly distributed
in $[0,10]$.
To each $P$ corresponds a polynomial $Q$, given by
\bdis
Q(z) = Iz^{5} + B_{4}z^{4} + \dots + B_{0} \; ,
\edis
where the $2m \times 2m$ matrices $B_{j}$, for $j=0,\dots ,4$, are defined as in Theorem~\ref{squarematrix}. 
For each value of $m$, one thousand such random matrix polynomials were generated. We then compared the ratios
of the bounds to the moduli of their largest and smallest eigenvalues, respectively.  
The following upper and lower bounds were compared:

\begin{tabular}{l l}
                              &                                                                                \\
          {\bf Upper bounds:} & $P$ :    Theorem~\ref{genCauchy} with $P$                               \\
                              & $Q$ :     Theorem~\ref{squarematrix} with $Q$ \\
                              &                                                                                \\
          {\bf Lower bounds:} & $P$ : Lower bound in Theorem~\ref{genCauchy} with $P$                \\
          & $Q$ : Lower bound in Theorem~\ref{squarematrix} with $Q$ \\
          & $A_{0}^{-1}P$ : Lower bound in Theorem~\ref{genCauchy} with $A_{0}^{-1}P$ instead of $P$                     \\
          & $B_{0}^{-1}Q$ : Lower bound in Theorem~\ref{squarematrix} with $B_{0}^{-1}Q$ instead of $Q$ \\
          & $Q_{R}$ : Reciprocal of upper bound in Theorem~\ref{squarematrix} with $Q_{R}$.
\end{tabular}

$ $ \newline Preceding each bound's description is its corresponding designation in the tables below. We recall that
the bound designated by $Q_{R}$ is obtained by forming $Q_{R}$ from $P_{r}$ just as $Q$ is formed from $P$.
This notation will also be used for all subsequent examples.

In Table~\ref{Cauchyupperbnds} we have listed as percentages the means of the ratios of the upper bound to 
the modulus of the largest eigenvalue with their standard deviation between parentheses. For instance, a 
ratio of $3/2$ corresponds to $150 \protect \%$. The closer the number is to $100$, the better it is, 
both for upper and lower bounds. The bounds were
computed with the $1$-Norm, $\infty$-Norm, and the $2$-Norm.
In Table~\ref{Cauchyupperfrqs} we have listed the number of times each bound was the better bound.
In Table~\ref{Cauchylowerbnds} and Table~\ref{Cauchylowerfrqs} the same was done for the lower bounds.

As Table~\ref{Cauchyupperbnds} shows, using Theorem~\ref{squarematrix} produces better upper bounds on average
for all three norms, except for scalar polynomials. This advantage seems to increase with increasing $m$.
The standard deviations are generally comparable, except for the higher values of $m$ with the $2$-Norm, where they are
larger for the $Q$ bound.
The 2-Norm gives the best results but is more costly to compute. It is also clear from Table~\ref{Cauchyupperfrqs}
that even when a bound is worse on average, there is still a non-negligible number of cases where that bound
prevails over the other one.

For the lower bounds, Table~\ref{Cauchylowerbnds} and Table~\ref{Cauchylowerfrqs} show that the advantage 
goes to $Q_{R}$ for $m > 1$, while $P$ dominates for $m=1$. As expected, the bounds $A_{0}^{-1}P$ and 
$B_{0}^{-1}Q$ are better than $P$ and $Q$, respectively. There were no instances where the latter two delivered 
the best bound.   
As for the upper bounds, the $2$-Norm gives the best results. 

Summarizing, we can say that, for $m \geq 2$, the upper bounds were, on average, improved by the use of 
Theorem~\ref{squarematrix}
with $Q$, whereas for the lower bounds, the same is true for Theorem~\ref{squarematrix} with $Q_{R}$ instead of $Q$.
The quality of the bounds improved with higher values of $m$. Although this is but one class of examples, the matrices
were randomly generated without any effort at special selection.


\begin{table}[H]
\begin{center}
\footnotesize
\begin{tabular}{c c c}
\begin{tabular}{c|c c}
 m                & $P$        & $Q$         \\  \hline 
1                 & 116 (17)   & 131 (18)    \\
2                 & 161 (32)   & 162 (29)    \\
10                & 316 (45)   & 247 (47)    \\
25                & 481 (62)   & 312 (67)    \\
\end{tabular}
&
\qquad
\begin{tabular}{c|c c}
 m                & $P$        & $Q$         \\  \hline 
1                 & 116 (17)   & 122 (18)    \\
2                 & 162 (33)   & 147 (30)    \\
10                & 317 (46)   & 216 (50)    \\
25                & 481 (62)   & 271 (71)    \\
\end{tabular}
&
\qquad
\begin{tabular}{c|c c}
 m                & $P$        & $Q$         \\  \hline 
1                 & 116 (17)   & 118 (16)    \\
2                 & 140 (28)   & 133 (23)    \\
10                & 167 (19)   & 161 (31)    \\
25                & 177 (14)   & 174 (39)    \\
\end{tabular}
\\
 & & \\
\fbox{1-Norm} & \fbox{$\infty$-Norm} & \fbox{2-Norm} 
\end{tabular}
\caption{Comparison of the upper bounds with $n=10$.}
\label{Cauchyupperbnds}
\end{center}
\end{table}
\normalsize

\vskip -0.5cm
\begin{table}[H]
\begin{center}
\footnotesize
\begin{tabular}{c c c}
\begin{tabular}{c|c c}
 m                & $P$        & $Q$         \\  \hline 
1                 & 939        & 61          \\
2                 & 583        & 417         \\
10                & 152        & 848         \\
25                & 93         & 907         \\
\end{tabular}
&
\qquad
\begin{tabular}{c|c c}
 m                & $P$        & $Q$         \\  \hline 
1                 & 759        & 241         \\
2                 & 275        & 725         \\
10                & 136        & 864         \\
25                & 90         & 910         \\
\end{tabular}
&
\qquad
\begin{tabular}{c|c c}
 m                & $P$        & $Q$         \\  \hline 
1                 & 557        & 443         \\
2                 & 342        & 658         \\
10                & 317        & 683         \\
25                & 340        & 660         \\
\end{tabular}
\\
 & & \\
\fbox{1-Norm} & \fbox{$\infty$-Norm} & \fbox{2-Norm} 
\end{tabular}
\caption{Best bound frequencies for the upper bounds with n=10.}
\label{Cauchyupperfrqs}
\end{center}
\end{table}
\normalsize

\begin{center}
\footnotesize
\begin{tabular}{c c}
\begin{tabular}{c|c c c c c}
 m               &  $P$       &  $Q$       &  $A_{0}^{-1}P$     &  $B_{0}^{-1}Q$    &   $Q_{R}$      \\  \hline 
1                &  87 (9)    &  33 (21)   &  87 (9)            &  43 (24)          &   81 (9)       \\
2                &  49 (13)   &  16 (13)   &  65 (13)           &  25 (18)          &   67 (11)      \\
10               &  10 (6)    &  4 (4)     &  27 (10)           &  8 (10)           &   39 (9)       \\
25               &  4 (3)     &  2 (2)     &  16 (7)            &  5 (6)            &   29 (8)       \\
\end{tabular}
&
\qquad
\begin{tabular}{c|c c c c c}
 m               &  $P$       &  $Q$       &  $A_{0}^{-1}P$     &  $B_{0}^{-1}Q$    &   $Q_{R}$      \\  \hline 
1                &  87 (9)    &  34 (21)   &  87 (9)            &  43 (24)          &   83 (9)       \\
2                &  49 (13)   &  16 (13)   &  65 (13)           &  26 (18)          &   70 (11)      \\
10               &  9 (6)     &  4 (5)     &  27 (10)           &  9 (10)           &   42 (10)      \\
25               &  4 (3)     &  2 (3)     &  16 (7)            &  5 (6)            &   31 (9)       \\
\end{tabular}
\\
&  \\
\fbox{1-Norm} & \fbox{$\infty$-Norm} 
\end{tabular}
\end{center}

\vskip -0.25cm
\begin{table}[H]
\begin{center}
\footnotesize
\begin{tabular}{c|c c c c c}
 m               &  $P$       &  $Q$       &  $A_{0}^{-1}P$     &  $B_{0}^{-1}Q$    &   $Q_{R}$      \\  \hline 
1                &  87 (9)    &  38 (23)   &  87 (9)            &  45 (25)          &   87 (8)       \\
2                &  61 (14)   &  21 (17)   &  71 (13)           &  29 (20)          &   78 (12)      \\
10               &  25 (10)   &  7 (9)     &  39 (12)           &  11 (13)          &   54 (12)      \\
25               &  15 (7)    &  4 (6)     &  26 (9)            &  6 (9)            &   42 (11)      \\
\end{tabular}
\\
\vskip 0.25cm
\fbox{2-Norm}  
\caption{Comparison of the lower bounds with $n=10$.}
\label{Cauchylowerbnds}
\end{center}
\end{table}

\vskip -0.5cm
\begin{table}[H]
\begin{center}
\footnotesize
\begin{tabular}{c c c}
\begin{tabular}{c|c c c}
 m           & $A_{0}^{-1}P$   & $B_{0}^{-1}Q$   &   $Q_{R}$       \\  \hline 
1            & 761             & 9               &  230            \\
2            & 395             & 19              &  586            \\
10           & 85              & 17              &  898            \\
25           & 48              & 11              &  941            \\
\end{tabular}
&
\qquad
\begin{tabular}{c|c c c}
 m           & $A_{0}^{-1}P$   & $B_{0}^{-1}Q$   &   $Q_{R}$       \\  \hline 
1            & 720             & 22              &  258            \\
2            & 277             & 19              &  704            \\
10           & 99              & 7               &  894            \\
25           & 56              & 3               &  941            \\
\end{tabular}
&
\qquad
\begin{tabular}{c|c c c}
 m           & $A_{0}^{-1}P$   & $B_{0}^{-1}Q$   &   $Q_{R}$       \\  \hline 
1            & 394             & 15              & 591             \\
2            & 175             & 11              & 814             \\
10           & 76              & 9               & 915             \\
25           & 61              & 3               & 936             \\
\end{tabular}
\\
 & & \\
\fbox{1-Norm} & \fbox{$\infty$-Norm} & \fbox{2-Norm} 
\end{tabular}
\caption{Best bound frequencies for the lower bounds with $n=10$.}
\label{Cauchylowerfrqs}
\end{center}
\end{table}
\normalsize

\vskip -1.0cm
$ $\newline \underline{\bf Example 2}
\newline \newline
In this example we compare the ability to produce a gap between the moduli of two groups of eigenvalues
for Theorem~\ref{genPellet} and Theorem~\ref{squarematrix}.
To do so we have generated a class of matrix polynomials of degree $n=14$ 
with $m \times m = 25 \times 25$ matrix coefficients of the form   
\bdis
P(z) = Iz^{14} + A_{13}z^{13} + \dots + A_{0} \; .
\edis
The elements of the matrices $A_{j}$ have real and imaginary parts that are uniformly distributed in   
$[-50^{2}/2,50^{2}/2]$, $[-200^{2}/2,200^{2}/2]$, and $[-2,2]$ for $j=11$, $j=12$, and $j\neq 11,12,13$, respectively.
One thousand such polynomials were generated for which the real and imaginary parts of the elements of $A_{13}$
are in the intervals $[-\eta,\eta]$, where $\eta=0, 0.25, 0.5, 1$. 
For each matrix polynomial, we verified if Theorem~\ref{genPellet} and Theorem~\ref{squarematrix} were applicable 
and if they were, computed the ratio of the gap between the moduli of the two groups of 
eigenvalues to the actual gap. The theorems were applied to both $P$ and $Q$, and 
also to both $A_{k}^{-1}P$ and $B_{k/2}^{-1}Q$, with $k=12$.

Throughout this example, only the $1$-Norm was used to limit the number of tables. Very similar results are 
obtained for the $\infty$-Norm. Better results for all bounds are achieved for the 2-Norm (spectral norm), 
but the relative performance of the bounds follows the same trends. The $2$-Norm's computational cost is higher than
the other two norms, but we do not advocate the use of any particular norm as this choice depends very 
much on the situation.

On the left in Table~\ref{ex2table1} are listed the number of times (out of $1000$ cases) 
that a gap was computed between the $mk=300$ smallest and $m(n-k)=50$ largest eigenvalues for each value of 
$\eta$ when applying Theorem~\ref{genPellet} and Theorem~\ref{squarematrix} to the matrix polynomials $P$ and $Q$. 
The table contains the total number 
of computed gaps for each matrix polynomial as well as the number of times each matrix polynomial was the 
only one of the two for which a gap could be computed. On the right in Table~\ref{ex2table1} are listed, as 
percentages, the means of the ratios of the gap from Theorem~\ref{genPellet} and Theorem~\ref{squarematrix}
to the actual gap, with the standard 
deviations between parentheses. The rightmost column lists the percentage of cases where the gap for 
$Q$ was larger than the one for $P$ when a gap was produced for both.
Table~\ref{ex2table2} is the analogous table obtained by applying the theorems to 
$A_{k}^{-1}P$ and $B_{k/2}^{-1}Q$ instead of $P$ and $Q$, respectively.

Overall, the results obtained by applying Theorem~\ref{squarematrix} improve as $||A_{13}||$ becomes smaller,
which is understandable since it tends to keep $||B_{j}||$ ($j=0,\dots ,n/2-1$) of the same order of magnitude
as $||A_{j}||$ ($j=0,\dots ,n-1$). This is more important here than for mere upper and lower bounds as in Example 1,
because the existence of positive roots for $f_{k}$ in Theorem~\ref{genPellet} and 
$g_{k}$ in Theorem~\ref{squarematrix} is heavily dependent on the relative magnitudes of their coefficients.
The results show clearly that, as $||A_{13}||$ becomes smaller, the application of Theorem~\ref{squarematrix}
improves the number of times a gap can be detected when compared to Theorem~\ref{genPellet} and 
there was a large number of cases in which a gap could be computed for $Q$, but not for $P$. 
The average gap width is generally better for $P$ and $A_{k}^{-1}P$ and becomes better for $Q$ only when $\eta \approx 0$.

For this example we found that when $\eta$ is small, the number of gaps detected for $Q$ 
is not much lower than for $B_{k/2}^{-1}Q$, which is more costly to compute.

%
%

%
%
\begin{table}
\begin{center}
\footnotesize       
\tabcolsep=0.50cm
\begin{tabular}{c c}
\tabcolsep=0.10cm
\begin{tabular}{c|cccc}
$\eta$
&\begin{tabular}{c} $P$  \\ TOTAL \end{tabular} 
& \begin{tabular}{c} $Q$  \\ TOTAL \end{tabular} 
& \begin{tabular}{c} $P$ \\ ONLY \end{tabular} 
& \begin{tabular}{c} $Q$ \\ ONLY \end{tabular} \\  \hline 
1    &     210              &        17                    &    193                   &       0                   \\
0.5  &     314              &        535                   &    0                     &       221                 \\
0.25 &     379              &        865                   &    0                     &       486                 \\
0    &     433              &        907                   &    0                     &       474                 \\
\end{tabular}   
&
\tabcolsep=0.20cm
\begin{tabular}{c|c c c}
$\eta$
& $P$  
& $Q$ 
& \begin{tabular}{c} \protect \%  \\  Gap(Q) > Gap(P) \end{tabular} \\ \hline 
1    &    16 (7)           &      5 (3)                   &        0                                \\
0.5  &    18 (8)           &      12 (5)                  &        22                               \\
0.25 &    21 (8)           &      19 (6)                  &        64                               \\
0    &    22 (8)           &      26 (8)                  &        100                              \\
\end{tabular}
\\
&  \\
\fbox{Gap frequency} & \fbox{Gap ratio}    
\end{tabular}
\caption{Gap frequencies and gap ratios for $P$ and $Q$ when $n=14$, $m=25$, and $k=12$.}    
\label{ex2table1}
\end{center}
\end{table}

\vskip -0.6cm

\begin{table}
\begin{center}
\footnotesize       
\tabcolsep=0.50cm
\begin{tabular}{c c}
\tabcolsep=0.10cm
\begin{tabular}{c|cccc}
$\eta$
&\begin{tabular}{c}  $A_{k}^{-1}P$    \\ TOTAL \end{tabular} 
& \begin{tabular}{c} $B_{k/2}^{-1}Q$     \\ TOTAL \end{tabular} 
& \begin{tabular}{c} $A_{k}^{-1}P$      \\ ONLY \end{tabular} 
& \begin{tabular}{c} $B_{k/2}^{-1}Q$    \\ ONLY \end{tabular} \\  \hline 
1    &     900              &       543                    &    357                   &       0                  \\
0.5  &     905              &       885                    &    20                    &       0                  \\
0.25 &     923              &       958                    &    1                     &       36                 \\
0    &     925              &       973                    &    0                     &       48                  \\
\end{tabular}   
&
\tabcolsep=0.20cm
\begin{tabular}{c|c c c}
$\eta$
& $A_{k}^{-1}P$  
& $B_{k/2}^{-1}Q$ 
& \begin{tabular}{c} \protect \%  \\  Gap($B_{k/2}^{-1}Q$) > Gap($A_{k}^{-1}P$) \end{tabular} \\ \hline 
1    &    30 (8)           &      9 (4)                   &        0                                \\
0.5  &    30 (8)           &      16 (5)                  &        1                                \\
0.25 &    31 (8)           &      23 (5)                  &        4                                \\
0    &    31 (8)           &      29 (7)                  &        17                               \\
\end{tabular}
\\
&  \\
\fbox{Gap frequency} & \fbox{Gap ratio}    
\end{tabular}
\caption{Gap frequencies and gap ratios for $A_{k}^{-1}P$ and $B_{k/2}^{-1}Q$ when $n=14$, $m=25$, and $k=12$.}    
\label{ex2table2}
\end{center}
\end{table}

\vskip -0.5cm
$ $ \newline \underline{\bf Example 3}
\newline \newline
This example illustrates how Theorem~\ref{squarematrix}, applied to scalar polynomials ($m=1$), can sometimes be used to
provide gaps between the moduli of zeros of a polynomial when the classical scalar version of Pellet's theorem cannot.
To this end, we have generated $1000$ scalar polynomials of the form
\bdis
p(z) = z^{20} + a_{19} z^{19} + \dots + a_{0} \; ,
\edis
where the coefficients $a_{j}$ are random real numbers, uniformly distributed on $[-2,-1] \cup [1,2]$ for $j=3,5,11,13$, 
on $[-10,-8] \cup [8,10]$ for $j=4$, on $[-16,-14] \cup [14,16]$ for $j=12$, and on 
$[-1,1]$ for all $j \neq 3,4,5,11,12,13$. 
Theorem~\ref{genPellet} was applied for $k=4,12$.
Because the matrices involved are all $2 \times 2$, we have used the 2-Norm, which can easily be computed in this case
and we have applied Theorem~\ref{squarematrix} only to $B_{k/2}^{-1}Q$ since here the matrix multiplications are not costly
and the results are better than when the theorem is applied to $Q$.

The results are displayed in Table~\ref{ex3table1} in the same format as for the previous example.
The designation $A_{k}^{-1}P$ is replaced by $p$ since for the scalar version of 
Pellet's theorem there is no difference between its application to $p$ and to $a_{k}^{-1}p$.
The numbers of times that a gap could be computed for both $k=4$ and $k=12$ for the same polynomial are listed 
in Table~\ref{ex3table2}.

For $k=4$, more gaps are detected by the classical Pellet theorem than with Theorem~\ref{squarematrix}, 
although there are still 21 cases where the latter does 
provide a gap when the classical Pellet theorem cannot and that gap is also wider in $29 \protect \%$ of the 
cases when both provide a gap. For $k=12$, the situation is reversed: using the generalized Pellet theorem 
for $B_{k/2}^{-1}Q$ delivers significantly better results than using the scalar Pellet theorem for $p$.
There seems to be no particular reason why the classical Pellet theorem is better for $k=4$, but not for $k=12$.

As Table~\ref{ex3table2} shows, $B_{k/2}^{-1}Q$ delivered a gap for both $k=4$ and $k=12$ in $30$ cases
more than the $4$ cases in which the classical Pellet theorem did.


\begin{table}
\begin{center}
\footnotesize
\tabcolsep=0.50cm
\begin{tabular}{c c}
\tabcolsep=0.10cm
\begin{tabular}{c|cccc}
$ k  $
&\begin{tabular}{c} $ p         $ \\ TOTAL \end{tabular} 
& \begin{tabular}{c}   $B_{k/2}^{-1}Q$ \\ TOTAL \end{tabular} 
& \begin{tabular}{c} $p$ \\ ONLY \end{tabular} 
& \begin{tabular}{c} $B_{k/2}^{-1}Q$ \\ ONLY \end{tabular} \\  \hline 
4    &     378              &        223                   &    176                   &       21                  \\
12   &     59               &        98                    &     9                    &       48                  \\
\end{tabular}   
&
\tabcolsep=0.20cm
\begin{tabular}{c|c c c}
$ k $
& $ p $                    
& $B_{k/2}^{-1}Q$ 
& \begin{tabular}{c} \protect \%  \\  Gap($B_{k/2}^{-1}Q$) > Gap($p$) \end{tabular} \\ \hline 
4    &    44 (18)          &      37 (15)                 &        29                               \\
12   &    26 (13)          &      28 (13)                 &        82                               \\
\end{tabular}
\\
&  \\
\fbox{Gap frequency} & \fbox{Gap ratio}    
\end{tabular}
\caption{Gap frequencies and gap ratios for $p$ and $B_{k/2}^{-1}Q$ when $n=20$, $m=1$, and $k=4,12$.}    
\label{ex3table1}      
\end{center}
\end{table}

\vskip -0.5cm
\begin{table}
\begin{center}
\footnotesize
\tabcolsep=0.50cm
\begin{tabular}{cccc}
\begin{tabular}{c} $ p         $ \\ TOTAL \end{tabular} 
 & \begin{tabular}{c}   $B_{k/2}^{-1}Q$ \\ TOTAL \end{tabular} 
 & \begin{tabular}{c} $p$ \\ ONLY \end{tabular} 
 & \begin{tabular}{c} $B_{k/2}^{-1}Q$ \\ ONLY \end{tabular} \\  \hline 
     23               &        49                    &     4                    &       30                  \\
\end{tabular}   
\caption{Gap frequencies when $n=20$ and $m=1$ for \emph{both} $k=4$ and $k=12$.}    
\label{ex3table2}      
\end{center}
\end{table}
\vskip -0.5cm 
$ $ \newline
{\bf Remark.} The reported average gap widths are over all cases that a gap was produced for 
a particular theorem and \emph{not} over only those cases where \emph{both} theorems produced a gap.
The difference in gap width is usually not by an order of magnitude, which would be unexpected.
The only time this happens is in limiting cases, where one theorem barely manages to
detect a gap and the other detects a gap of average width. 
Our interest was less in the gap width than
in detecting gaps that could previously not be detected.

%
%
%
%
\section{Scalar polynomials as matrix polynomials}
\label{scalarpol}      

In the remainder of this work, we suggest an idea for an alternative way to treat scalar polynomials.
More specifically, we propose expressing a scalar polynomial as the determinant of a matrix polynomial 
to which the generalized Pellet theorem can then be applied, instead of applying the regular (scalar) 
version of Pellet's theorem to the original scalar polynomial. The intention is, as before, to try and 
overcome situations where the real polynomial $f_{k}$ in Theorem~\ref{Pellet} does not have positive roots,
and to create new upper and lower bounds on the moduli of the zeros with the generalized Cauchy result
(Theorem~\ref{genCauchy}).

There are infinitely many ways to write a scalar polynomial as the determinant of a matrix polynomial 
and the best way to proceed is, a priori, not clear. As an example, consider the polynomial $z^4-z^2 +3z -2$. 
It can be written as
\bdis
\det
\begin{pmatrix}
z^{2} -\sqrt{2}  & z         \\
z-3      & z^{2} + \sqrt{2}  \\ 
\end{pmatrix} 
=
\det \left (
\begin{pmatrix}
1 & 0 \\
0 & 1 \\
\end{pmatrix} 
z^{2} + 
\begin{pmatrix}
0 & 1 \\
1 & 0 \\
\end{pmatrix} 
z + 
\begin{pmatrix}
-\sqrt{2}  & 0         \\
-3      &  \sqrt{2}  \\ 
\end{pmatrix} 
\right )
\; ,
\edis
but it can also be written as, e.g.,  
\bdis
\det
\begin{pmatrix}
z^{2}   & -1        \\
3z-2     & z^{2} -1          \\ 
\end{pmatrix} 
=
\det \left (
\begin{pmatrix}
1 & 0 \\
0 & 1 \\
\end{pmatrix} 
z^{2} + 
\begin{pmatrix}
0 & 0 \\
3 & 0 \\
\end{pmatrix} 
z + 
\begin{pmatrix}
0 & -1  \\                                 
-2  &  -1  \\ 
\end{pmatrix} 
\right )
\; .
\edis
Furthermore, one is not limited to using $2 \times 2$ matrices, and matrices of larger size can also be used.
Although this approach is not applicable for every value of $k$ in Theorem~\ref{Pellet}, it does
provide new upper and/or lower bounds on the moduli of the zeros.

It is difficult to predict which of the equivalent matrix polynomial representations of a general 
scalar polynomial will yield the best results, but numerical experiments seem to indicate that better 
results are obtained if the coefficients are, in a sense, equally distributed among the elements
of the matrix polynomial. This distribution therefore needs to be tailored to the particular polynomial 
under consideration, which makes it difficult to formulate a general method. Consequently, it is a conjecture 
that our idea might lead to a useful general method.

To strengthen this conjecture, we present the following lemma, where we apply it to           
a class of lacunary polynomials, followed by an illustrative numerical example.

\begin{lemma}
\label{lacpol} 
Let
\beq
\label{lacpoldef1}           
p(z) = az^{n}+bz^{n-1}+cz^{n-2}+\alpha z^{2}+\beta z+\gamma \; ,
\eeq
where $a,b,c,\alpha,\beta,\gamma \in \complex$ with $a\alpha \neq 0$, and let
\bdis
A =
\begin{pmatrix}
\sqrt{a} & 0         \\
0        & \sqrt{a}  \\ 
\end{pmatrix} 
\;\; ,  \;\;
B = \begin{pmatrix}
\dfrac{b}{2\sqrt{a}}+ i \lb c - \dfrac{b^{2}}{4a} \rb^{1/2}   &      0                        \\
0                        &  \dfrac{b}{2\sqrt{a}}- i \lb c - \dfrac{b^{2}}{4a} \rb^{1/2}      \\ 
\end{pmatrix}
\;\; ,  \;\;
\edis
\bdis
C =
\begin{pmatrix}
\sqrt{a} & 0         \\
0        & 0         \\ 
\end{pmatrix} 
\;\; ,  \;\;
D =
\begin{pmatrix}
\dfrac{b}{2\sqrt{a}}+ i \lb c - \dfrac{b^{2}}{4a} \rb^{1/2}        & 0         \\
0        & \sqrt{a}  \\ 
\end{pmatrix} 
\;\; ,  \;\;     
E = 
\begin{pmatrix}
0              &   0                \\
0              & \dfrac{b}{2\sqrt{a}}- i \lb c - \dfrac{b^{2}}{4a} \rb^{1/2}      \\
\end{pmatrix}
\; ,
\edis
\bdis
V = 
\begin{pmatrix}
0              &   -\sqrt{\alpha}   \\
\sqrt{\alpha}  &    0  \\
\end{pmatrix}
\;\; , \; \text{and} \;\; 
W = 
\begin{pmatrix}
0              &   -\dfrac{\beta}{2\sqrt{\alpha}}- i \lb \gamma - \dfrac{\beta^{2}}{4\alpha} \rb^{1/2}    \\
\dfrac{\beta}{2\sqrt{\alpha}}- i \lb \gamma - \dfrac{\beta^{2}}{4\alpha} \rb^{1/2}    &    0  \\
\end{pmatrix}
\;\; .
\edis
Then the zeros of $p$ are the eigenvalues of the matrix polynomial 
\beq
\label{lacpoldef2}           
Q_{even}(z) = A z^{n/2} + B z^{n/2-1} + V z + W  
\eeq
when its degree $n$ is even, and they are the finite eigenvalues of
\beq
\label{lacpoldef3}           
Q_{odd}(z) =  C z^{(n+1)/2} + D z^{(n-1)/2} + E z^{(n-3)/2} + V z + W  
\nonumber
\eeq
when its degree $n$ is odd.
\end{lemma}

\prf 
The polynomial $p$ can be written as follows:
\begin{eqnarray*}
p(z) & = & az^{n}+bz^{n-1}+cz^{n-2}+\alpha z^{2}+\beta z+\gamma  \\
     & = & az^{n-2}\lb z^{2}+\dfrac{b}{a} z + \dfrac{c}{a} \rb 
           + \alpha \lb z^{2}+\dfrac{\beta}{\alpha}z + \dfrac{\gamma}{\alpha} \rb \\ 
     & = & z^{n-2}\lb a \lb z+\dfrac{b}{2a} \rb^{2}  + c - \dfrac{b^{2}}{4a} \rb 
           + \lb \alpha \lb z+\dfrac{\beta}{2\alpha} \rb^{2} + \gamma - \dfrac{\beta^{2}}{4\alpha} \rb \; . \\ 
\end{eqnarray*}
Now assume that $n$ is even. After factoring the quadratic expressions, we obtain
\beq
\label{bigdeteven}
p(z) = \text{det}
\begin{pmatrix}
z^{n/2-1}\lb \sqrt{a}z+\dfrac{b}{2\sqrt{a}}+ i \lb c - \dfrac{b^{2}}{4a} \rb^{1/2}  \rb & 
- \sqrt{\alpha}z-\dfrac{\beta}{2\sqrt{\alpha}}- i \lb \gamma - \dfrac{\beta^{2}}{4\alpha} \rb^{1/2}   \\ 
\sqrt{\alpha}z+\dfrac{\beta}{2\sqrt{\alpha}}- i \lb \gamma - \dfrac{\beta^{2}}{4\alpha} \rb^{1/2}  &   
z^{n/2-1}\lb \sqrt{a}z+\dfrac{b}{2\sqrt{a}}- i \lb c - \dfrac{b^{2}}{4a} \rb^{1/2}  \rb \\
\end{pmatrix}
\; ,
\nonumber
\eeq
which is easily seen to be equivalent to $\text{det}(Q_{even}(z))$, with   
\bdis
Q_{even}(z) = A z^{n/2} + B z^{n/2-1} + V z + W   \; ,
\edis
where the $2 \times 2$ matrices $A$, $B$, $V$, and $W$ are defined in the statement of the lemma.

When $n$ is odd, we can similarly write $p$ as
\beq
\label{bigdetodd}   
p(z) = \text{det}
\begin{pmatrix}
z^{(n-1)/2}\lb \sqrt{a}z+\dfrac{b}{2\sqrt{a}}+ i \lb c - \dfrac{b^{2}}{4a} \rb^{1/2}  \rb & 
- \sqrt{\alpha}z-\dfrac{\beta}{2\sqrt{\alpha}}- i \lb \gamma - \dfrac{\beta^{2}}{4\alpha} \rb^{1/2}   \\ 
\sqrt{\alpha}z+\dfrac{\beta}{2\sqrt{\alpha}}- i \lb \gamma - \dfrac{\beta^{2}}{4\alpha} \rb^{1/2}  &   
z^{(n-3)/2}\lb \sqrt{a}z+\dfrac{b}{2\sqrt{a}}- i \lb c - \dfrac{b^{2}}{4a} \rb^{1/2}  \rb \\
\end{pmatrix}
\; ,
\nonumber
\eeq
which is equivalent to $\text{det}(Q_{odd}(z))$, with    
\bdis
Q_{odd}(z) = C z^{(n+1)/2} + D z^{(n-1)/2} + E z^{(n-3)/2} + V z + W   \; ,
\edis
where the $2 \times 2$ matrices $C$, $D$, $E$, $V$, and $W$ 
are defined in the statement of the lemma.
\qed

Lemma~\ref{lacpol} can now be used in conjunction with Theorem~\ref{genPellet} and Theorem~\ref{genCauchy} 
to obtain additional and/or improved results for the scalar versions of these theorems, as long as the 
appropriate matrices are nonsingular. In this regard we note that the matrices $A$ and $V$ are always 
nonsingular since we assumed that $a \alpha \neq 0$, and the matrices $C$ and $E$ are always singular.
As with Theorem~\ref{squarematrix}, not all of the scalar bounds can be improved, but enough can to make it worthwile.
When applying a similar process as in Lemma~\ref{lacpol} to a general polynomial, 
the bounds, for a particular value of $k$, require ${\mathcal{O}}(n)$ operations since the coefficient matrices 
of the equivalent matrix polynomial are always of size $2 \times 2$.
Similarly to what we had before, the degrees of the new matrix polynomials $Q_{even}$ and $Q_{odd}$ are half
and roughly half that of $P$, respectively.
$ $ \newline
{\bf Remarks.} 
Although we will not pursue the matter further, Lemma~\ref{lacpol} allows for many variations on its theme. 
For instance, a polynomial of the form $az^{n}+bz^{n-1}+cz^{n-2}+\alpha z + \beta$ with $a\alpha \neq 0$ becomes 
of the form of a polynomial as in Lemma~\ref{lacpol} with $\gamma=0$ after multiplication by $z$, so that the 
lemma becomes applicable to such polynomials as well.
Moreover, the lacunary polynomial $p$ in Lemma~\ref{lacpol} is a special case for $r=s=1$ of the polynomial
\bdis
az^{n+2r-2}+bz^{n+r-2}+cz^{n-2}+\alpha z^{2s}+\beta z^{s}+\gamma \; ,
\edis
where $r \geq 1$ and $s \geq 1$ are integers. Lemma~\ref{lacpol} can easily be modified to include these more 
general polynomials, which can be further generalized by replacing $z$ with a polynomial in $z$.

The following example illustrates how Lemma~\ref{lacpol} can be used to improve upper and lower bounds
on the moduli of the zeros of a polynomial and also to detect gaps between groups of zeros where the classical
scalar Pellet theorem cannot, or to improve the width of the gaps that it does detect. 
$ $ \newline \newline \underline{\bf Example 4}
\newline \newline
For this example, we have generated
$1000$ polynomials of the form~(\ref{lacpoldef1}), with real random coefficients $a$, $b$, $c$, $\alpha$, $\beta$, and $\gamma$ 
uniformly distributed in $[-50,50]$ and for $n=20,40,80$. Both the upper and lower bounds from Theorem~\ref{genCauchy} 
and its scalar version were computed, as well as the bounds from Theorem~\ref{genPellet} and its scalar version
for $k=2$ and $k=n-2$, i.e., bounds on the gap between $2$ and $n-2$ zeros.
These bounds were computed for both $p$ and $Q_{even}$. The matrices $B$ and $W$ defining $Q_{even}$ 
were nonsingular in all cases.
Again, as before, because all matrices involved are $2 \times 2$, we used only the $2$-Norm.
We have computed means and standard deviations of ratios and expressed them as percentages and we have listed gap
frequencies with the same conventions as in the previous examples. The results were collected in the tables below.  

Table~\ref{ex6table1} lists the ratios of the upper and lower bounds to the moduli of the largest and smallest zeros, 
respectively. The designation "scalar" and "matrix" refers to the use of $p$ from~(\ref{lacpoldef1}) and 
$Q_{even}$ from~(\ref{lacpoldef2}), respectively.
The last three columns list the percentage of cases for which using $Q_{even}$ yielded a better result than using $p$
for the upper, lower and both upper and lower bounds on the moduli of the zeros of $p$.
The upper and lower bound ratios using $p$ are quite unaffected by the degree of $p$, whereas using $Q_{even}$
yields better ratios and smaller standard deviations that improve with increasing degree $n$, a property that 
seems to be independent of the distribution of the coefficients.

Table~\ref{ex6table2} lists the number of times a gap could be computed with both the scalar and generalized (matrix) 
version of Theorem~\ref{genPellet} for $k=2$ and $k=n-2$. There is a clear advantage to using $Q_{even}$ that 
becomes more significant with increasing degree $n$. 
The same is true for the number of times both gaps could be computed, which can be found in Table~\ref{ex6table3}. To 
pick but one instance, for $n=40$, the classical scalar version of Pellet's theorem detects a gap
between the $38$ smallest and $2$ largest zeros $102$ times. The generalized Pellet theorem with $Q_{even}$
adds another $55$ times to that and also detects more than twice as many gaps (32 vs. 14) for both $k=2$ and $k=38$. 

Table~\ref{ex6table4} contains the means and standard deviations of the gaps and the percentage of times $Q_{even}$
produced a larger gap than $p$ for those cases in which both produced a gap. 
This was close to half of all cases, while producing a larger average ratio of computed gap 
to exact gap.
The computed gap may be larger for one or the other, although the situation could be reversed for
the bounds defining that gap. In other words, $Q_{even}$ may produce a smaller gap than $p$,
but its lower gap bound or its upper gap bound may be better. This means that combining the bounds from $Q_{even}$
and $p$ would improve the gap ratio even more.  


\begin{table}
\begin{center}
\footnotesize
\begin{tabular}{c||cc||cc||ccc}
$n$
& \begin{tabular}{c} Upper  \\ scalar  \end{tabular} 
& \begin{tabular}{c} Upper  \\ matrix  \end{tabular} 
& \begin{tabular}{c} Lower  \\ scalar  \end{tabular} 
& \begin{tabular}{c} Lower  \\ matrix  \end{tabular}  
& \begin{tabular}{c} $\protect \%$ Upper \\ better  \end{tabular}  
& \begin{tabular}{c} $\protect \%$ Lower \\ better  \end{tabular} 
& \begin{tabular}{c} $\protect \%$ Upper and Lower \\ better  \end{tabular} \\  \hline 
20     & 118 (28) & 108 (7) & 90 (15)  & 93 (7)  &  46     &   42    &   20         \\
40     & 119 (29) & 103 (4) & 90 (16)  & 98 (4)  &  65     &   59    &   39         \\
80     & 118 (31) & 102 (2) & 90 (16)  & 99 (3)  &  72     &   71    &   51         \\
\end{tabular}   
\caption{Upper and lower bounds for $n=20,40,80$.}                                                             
\label{ex6table1}
\end{center}
\end{table}
\vskip -0.5cm
\begin{table}
\begin{center}
\footnotesize
\begin{tabular}{c||cccc||cccc}
$n$
& \begin{tabular}{c} $k=2$   \\ scalar  \end{tabular} 
& \begin{tabular}{c} $k=2$   \\ matrix  \end{tabular} 
& \begin{tabular}{c} $k=2$   \\ scalar \\ only \end{tabular} 
& \begin{tabular}{c} $k=2$   \\ matrix \\ only \end{tabular} 
& \begin{tabular}{c} $k=n-2$ \\ scalar  \end{tabular} 
& \begin{tabular}{c} $k=n-2$ \\ matrix  \end{tabular}  
& \begin{tabular}{c} $k=n-2$ \\ scalar \\ only \end{tabular} 
& \begin{tabular}{c} $k=n-2$ \\ matrix \\ only \end{tabular} \\ \hline 
20     & 81   & 73   &  30  & 22   & 66   & 61   & 31  &  26   \\
40     & 123  & 161  & 19   & 57   & 102  & 157  & 18  &  73   \\
80     & 127  & 192  & 15   & 80   & 147  & 221  & 11  &  85   \\
\end{tabular}   
\caption{Gap frequency for $n=20,40,80$ and $k=2,n-2$.}
\label{ex6table2}
\end{center}
\end{table}
\vskip -0.5cm
\begin{table}
\begin{center}
\footnotesize
\begin{tabular}{c||cccc||cccc}
$n$
& \begin{tabular}{c} $k=2$ and $k=n-2$   \\ scalar  \end{tabular} 
& \begin{tabular}{c} $k=2$ and $k=n-2$  \\ matrix  \end{tabular} 
& \begin{tabular}{c} $k=2$ and $k=n-2$  \\ scalar only \end{tabular} 
& \begin{tabular}{c} $k=2$ and $k=n-2$  \\ matrix only \end{tabular} \\ \hline 
20     &  8   &  7   &  6   &   5     \\
40     & 14   & 32   &  4   &   22    \\
80     & 26   & 45   &  4   &   23    \\
\end{tabular}   
\caption{Gap frequency for $n=20,40,80$ and \emph{both} $k=2$ and $k=n-2$.}
\label{ex6table3}
\end{center}
\end{table}
\vskip -0.5cm
\begin{table}
\begin{center}
\footnotesize
\begin{tabular}{c||ccc||ccc}
$n$
& \begin{tabular}{c} $k=2$   \\ scalar  \end{tabular} 
& \begin{tabular}{c} $k=2$   \\ matrix  \end{tabular} 
& \begin{tabular}{c} $k=2$   \\ $\protect \%$ matrix better \\  \end{tabular} 
& \begin{tabular}{c} $k=n-2$ \\ scalar  \end{tabular} 
& \begin{tabular}{c} $k=n-2$ \\ matrix  \end{tabular}  
& \begin{tabular}{c} $k=n-2$ \\ $\protect \%$ matrix better \end{tabular} \\ \hline 
20     & 70 (23) & 66 (19)  &  44      &  68 (28)  &  71 (19) &  46          \\
40     & 74 (26) & 75 (21)  &  47      &  72 (29)  &  80 (19) &  48          \\
80     & 77 (29) & 85 (15)  &  46      &  71 (31)  &  85 (17) &  54          \\
\end{tabular}   
\caption{Gap ratio for $n=20,40,80$ and $k=2,n-2$.}                                                             
\label{ex6table4}
\end{center}
\end{table}

Similar results are obtained for odd powers and it was also observed that Lemma~\ref{lacpol} 
produces better results for this class of polynomials than Theorem~\ref{squarematrix}.


\end{document}